\documentclass[a4paper,10pt]{article}
\usepackage[frenchb]{babel}
\usepackage{amsfonts}
\usepackage{amssymb}
\usepackage{amsmath}
\usepackage{cite}
\newtheorem{theo}{Th\'eor\`eme}[section]
\newtheorem{lema}[theo]{Lemma}

\newtheorem{prop}[theo]{Proposition}

\newtheorem{coro}[theo]{Corolaire}
	
\newcommand{\gaminf}{\Gamma^\infty(\T,\OM)}
\newcommand{\Cinf}{\cinf(\T,\C)}

\newcommand{\LL}{\mathfrak{L}}

\newcommand{\M}{\mathbb{M}}

\newcommand{\Q}{\mathbb{Q}}
\newcommand{\R}{\mathbb{R}}
\newcommand{\D}{\mathbb{D}}

\newcommand{\T}{\mathbb{T}^1}
\newcommand{\cinf}{C^\infty}
\newcommand{\OM}{\Omega(\overline{\D})}
\newcommand{\C}{\mathbb{C}}
\newcommand{\N}{\mathbb{N}}
\newcommand{\eps}{\varepsilon}

\newcommand{\tal}{\theta+\alpha}

\newcommand{\Z}{\mathbb{Z}}

\newcommand{\Pru}{\mathbf{Preuve. }}
\newcommand{\CD}{\mathbb{CD}}

\begin{document}
\title{Sur la persistance des courbes invariantes pour les dynamiques holomorphes fibr\'ees lisses}
\author{Mario PONCE\\ Pontificia Universidad Cat\'{o}lica de Chile}
\maketitle
\begin{abstract}
En s'appuyant sur un th\'eor\`eme  des fonctions implicites de Hamilton nous montrons la persistance d'une courbe invariante indiff\'erente pour une dynamique holomorphe fibr\'ee en classe $\cinf$. Une condition diophantienne sur la paire de nombres de rotation est demand\'ee. On montre aussi que cette condition est optimale.
\end{abstract}
\section{Introduction}
La th\'eorie des syst\`emes dynamiques unidimensionnels n'est pas du tout un champ epuiss\'e. Cependant, il y en a au moins deux types des tels dynamiques qu'on pourrait consid\'erer comme \'etant bien comprises. Il s'agit de la th\'eorie des hom\'eomorphismes positifs du cercle (une dimension r\'eelle, voir par example\cite{HERM79},\cite{YOCC99}) et la dynamique semi-locale des germes holomorphes irrationnellement indiff\'erents qui fixent l'origine du plan complexe (une dimension complexe, voir par example \cite{CAGA93}, \cite{YOCC95}). On dit que telles dynamiques sont de type \emph{elliptique} car le comportement des orbites est contr\^ol\'e par un nombre r\'eel qu'on appele \emph{le nombre de rotation}. On sait que lorsque ce nombre de rotation v\'erifie une certaine  hypoth\`ese arithm\'etique alors la dynamique est plut\^ot simple (peu de mesures invariantes, conjugu\'ee \`a une dynamique lin\'eaire, etc.). Par contre, lorsqu'il s'agit de dynamiques d'allure elliptique en dimension sup\'erieure la th\'eorie n'est pas dans le m\^eme \'etat de d\'eveloppement. Dans  ce travail nous pr\'esentons des r\'esultats pour une classe de syst\`emes dynamiques qu'ont une allure elliptique, qui ne sont pas unidimensionnells et qu'on peut les consid\'erer comme \'etant de dimension interm\'ediaire. Ces objets sont une g\'en\'eralisation non triviale des germes holomorphes. D\'efinissons alors nos objets d'\'etude:  soit $\alpha$ un nombre irrationnel, qui nous fixons pour tout le reste de ce travail et  soit $U\subset \C$ un ouvert simplement connexe du plan complexe. Nous consid\'erons des transformations fibr\'ees de la forme
\begin{eqnarray*} 
F:\T\times U &\longrightarrow& \T\times \C\\
 (\theta,z)&\longmapsto& \big( \tal,f_{\theta}(z)\big),
 \end{eqnarray*}
 o\`u $\T$ est le cercle $\R/\Z$ et pour tout $\theta\in \T$ la fonction $f(\theta,\cdot):U\to\C$ est une fonction univalente (holomorphe et injective). Nous appelons  \`a  une telle transformation $F$ une \emph{dynamique holomorphe fibr\'ee au dessus du cercle} et la notons d\'esormais par \emph{dhf}. Dans ce travail nous allons traiter avec des \emph{dhf} qui sont de classe $\cinf$, dites \emph{dhf lisses}. 
 
 Ce type de transformations s'inscrit dans le  cadre plus g\'en\'eral des \emph{produits crois\'{e}s} (skew-products). Les produits crois\'es engendr\'es par des dynamiques qui sont bien comprises ont attir\'e  un grand int\'er\^et pendant les derni\`eres ann\'ees  car ils pr\'esentent parfois des ph\'enom\`enes nouveaux (voir \cite{HERM83}, \cite{VIAN97}, etc.). Dans le cas des deux types de dynamiques elliptiques unidimensionelles cit\'ees ci-dessus, la th\'eorie la mieux d\'evelopp\'ee est celle des ho\-m\'{e}omor\-phis\-mes du cercle fibr\'es au dessus d'une rotation irrationnelle.  L'article de M. Herman \cite{HERM83} pose les bases  de l'\'etude des hom\'{e}omorphismes fibr\'es en d\'efinisant le nombre de rotation fibr\'e.  Plus r\'{e}cemment cette  \'etude a \'et\'e relanc\'{e}e notamment par les travaux de G.Keller, T. J\"ager et J. Stark (voir \cite{STAR03},\cite{JAKE06}), qui ont \'etabli en particulier une classification \`a la Poincar\'e de telles dynamiques bas\'ee sur les propri\'et\'es des nombres de rotation associ\'es et l'existence de courbes invariantes.  Dans sa th\`ese de doctorat  O. Sester \cite{SEST97} (voir aussi \cite{SEST99})  a \'etudi\'e la dynamique des polyn\^omes fibr\'es, en g\'{e}n\'{e}ralisant les notions classiques d'ensemble de Julia, fonction de Green et de la cardio\"{i}de principale de l'ensemble de Mandelbrot dans l'espace de param\`etres. Les travaux de M. Jonsson (voir \cite{JONS99},\cite{JONS00}) sont aussi une r\'eference importante sur la dynamique fibr\'ee des transformations rationnelles. La reducibilit\'e des cocycles quasip\'eriodiques est un important sujet d'\'etude qui doit \^etre consid\'er\'e au moment de regarder les dynamiques fibr\'ees au dessus d'une rotation irrationnelle (voir par exemple \cite{KRIK99}, \cite{AVKR06}).  En fait, R.Johnson and J.Moser \cite{JOMO82} ont d\'efini un nombre de rotation fibr\'e pour l'op\'erateur de Schr\"odinger quasip\'eriodique. 
 
 Il y a encore une autre raison pour laquelle les \emph{dhf} pr\'esentent un int\'er\^et d'\'etude. Le probl\`eme qui nous traitons dans ce travail pr\'esente des characteristiques tr\`es similaires \`a celles de l'un des plus importants probl\`emes de la th\'eorie KAM, le bien connu probl\`eme de Melnikov (voir \cite{MELN65}, \cite{ELIA88}, \cite{BOUR97}). En fait, la persistance des tores invariants de dimension non-maximale pour les dynamiques Hamiltoniennes pr\'esente au niveau infinitesimal (et \emph{a posteriori} au niveau arithm\'etique) les m\^emes ingredients que la persistance des courbes invariantes pour les \emph{dhf}.  
 \\
  
  La notion de point fixe ou de point p\'eriodique pour une transformation fibr\'ee au dessus d'une rotation minimale n'a aucun sens (car $\alpha$ est un nombre irrationnel). L'extension naturelle de ce concept dans notre cadre est celle d'une courbe $u:\T\to \C$ invariante, c'est \`a dire une courbe de classe $\cinf$ qui satisfait l'\'equation 
 \begin{equation*}
 F\big(\theta,u(\theta)\big)=\big(\theta+\alpha,u(\theta+\alpha)\big)
 \end{equation*}
 pour tout $\theta\in \T$, ou de fa\c con equivalante $f\big(\theta,u(\theta)\big)=u(\tal)$. Ces objets jouent le r\^ole d'un centre autour duquel la dynamique de $F$ s'organise, en g\'en\'eralisant le r\^ole d'un point fixe pour la dynamique locale d'un germe holomorphe $g:(\C,0)\to (\C,0)$ (voir  \cite{PONC07-1}).
 
 Nous nous int\'eressons en l'existence d'une telle courbe invariante. Comme le nombre de rotation transversal (cf. d\'efinition \ref{deffrn}) contr\^ole la dynamique locale, nous nous concentrons en l'existence des courbes invariantes indiff\'erentes avec un nombre de rotation transversal fix\'e. Dans ce travail nous traitons le probl\`eme de la persistance des courbes, sous des petites perturbations sur la dynamique. Nous obtenons un r\'esultat du type KAM qui montre que la persistance a lieu sous certaines hypoth\`eses arithm\'etiques sur la paire des nombres de rotation associ\'ee. Ainsi, le probl\`eme est un probl\`eme de petits divis\'eurs. Dans la classe analytique nous pouvons obtenir des r\'esultats beaucoup plus int\'eressants. En fait, dans \cite{PONC07-3} l'auteur montre que la persistance a lieu sous une l'hypoth\`ese plus faible \`a la Brjuno sur la paire de nombres de rotation.
 
  \paragraph{Remerciements. } Ce travail fait partie de la Th\`ese de Doctorat de l'auteur, prepar\'ee au sein du Laboratoire de Math\'ematiques d'Orsay. L'auteur voudrait remercier  Jean-Christophe Yoccoz par sa direction et soutien constant. L'auteur remercie aussi Rapha\"el Krikorian pour ses corrections et suggestions. Ce travail a \'et\'e financ\'e par la Beca Master-Investigaci\'on y Doctorado CONICYT-Chile--Ambassade de France au Chili. Cet article a \'et\'e prepar\'e pendant un s\'ejour post-doctoral \`a l'Universidad Cat\'olica de Chile, suport\'e par le Proyecto ADI 17 Anillo en Sistemas Din\'amicos en baja dimensi\'on.  
  \section{Conditions arithm\'etiques et l'\'equation lin\'eairis\'ee}\label{sec2}
   Dans cette section nous allons d'abord rappeler quelques d\'efinitions classiques sur les conditions arithm\'etiques. Nous allons introduire  des conditions arithm\'etiques adapt\'ees aux probl\`emes de ce travail et nous traiterons aussi des \'equations lin\'eaires dites \emph{cohomologiques} qui font appel aux conditions arithm\'etiques.  
    \subsection{Conditions arithm\'etiques}\label{APPParitmetico}
Soit $x$ un nombre r\'eel. On pose
\begin{equation*}
\|x\|=\min_{p\in \Z} |x-p|
\end{equation*}
la distance au plus proche entier (ou bien la distance \`a l'origine $0$ dans $\T$ mesur\'ee sur le cercle). 
\paragraph{Conditions arithm\'etiques sur un nombre r\'eel. } Un traitement plus complet peut se trouver dans  plusieurs textes, dont \cite{LANG66}, \cite{CASS57}, \cite{SCHM80}. 
Soient $c>0$, $\tau\geq 0$. Nous d\'{e}finissons les ensembles de nombres r\'eels suivants 
\begin{displaymath}
		 \begin{array}{lll}
	\CD(c,\tau)&=&\big\{\alpha\in \T\setminus \Q\   \big|\  \forall N\geq 1\quad \|N\alpha\|\geq \frac{c}{N^{1+\tau}}\big\}\\
\\
	\CD(\tau)&=&\bigcup_{c>0}\CD(c,\tau)\\
\\
	\CD&=&\bigcup_{\tau\geq 0}\CD (\tau).\\

			\end{array}
\end{displaymath}
Si $\alpha$ appartient \`a $\CD$ on dit que $\alpha$ v\'erifie une condition diophantienne. Pour tout $\tau>0$ l'ensemble $\CD(\tau)$ est de mesure pleine au sens de Lebesgue. D'autre part, les ensembles $\CD(\tau)$ peuvent s'\'ecrire comme une union d\'enombrable d'ensembles ferm\'{e}s et d'int\'erieur vide, ce qui implique que ces ensembles sont petits du point de vue de la topologie (cat\'{e}gorie de Baire).
\paragraph{Conditions arithm\'etiques sur une paire}
Soient $\alpha, \beta$ des nombres r\'eels. Soient $c>0, \tau\geq 0$. Nous d\'{e}finissons les ensembles suivants
\begin{displaymath}
\begin{array}{lll}
\CD_1(c,\tau)&=&\big\{(\alpha,\beta)\in \T\times \T\   \big|\   \alpha \in \CD \textrm{ et }\forall N\in \Z \quad \|N\alpha-\beta\|\geq \frac{c}{N^{1+\tau}} \big\}\\
\\
\CD_1(\tau)&=&\bigcup_{c>0}\CD_1(c,\tau)\\
\\
\CD_1&=&\bigcup_{\tau\geq 0} \CD_1(\tau).\\
\end{array}
\end{displaymath}
Soit $\alpha\in \CD$.  L'ensemble  $\CD_1^{\alpha}(c,\tau)$ est l'nsemble de tous les $\beta$ tels que la paire $(\alpha,\beta)$ appartient \`a $\CD_1(c,\tau)$. De fa\c con analogue on definit $\CD_1^{\alpha}(\tau), \CD^{\alpha}_{1}$. Pour tout $\tau>0$ l'ensemble $\CD_1^{\alpha}(\tau)$ est de mesure pleine. De plus, l'ensemble $\CD_1^{\alpha}(c,\tau)$ est un ferm\'e d'int\'erieur vide.
Les ensembles $\CD_1^{\alpha}(\tau), \CD_1^{\alpha}$ s'\'ecrivent comme r\'eunions d\'enombrables d'ensembles du type $\CD_1^{\alpha}(c,\tau)$, et donc ils sont petits au sens topologique, en particulier leur compl\'{e}mentaire est dense dans le cercle. 
\subsection{\'{E}quation lin\'eairis\'ee}\label{APPPcohomologique}
Dans cette section nous allons rappeler des faits  bien connus sur l'\'equation cohomologique, nous renvoyons le lecteur aux textes \cite{HERM79}, \cite{KATO95}. 
Soit $\alpha$ dans $\R\setminus \Q$. Soit $\phi:\T\to \C$ une fonction de classe $\cinf$ de moyenne $\int_{\T}\phi(\theta)d\theta$ nulle pour la mesure de Lebesgue $d\theta$ du cercle. Nous cherchons une solution  $\psi:\T\to \C$ de classe $\cinf$ \`a l'\'equation \emph{cohomologique classique}
\begin{equation}\label{appppeqco}
\psi(\tal)-\psi(\theta)=\phi(\theta).
\end{equation}
Notons que la condition sur la moyenne de $\phi$ est n\'ecessaire. Une \'eventuelle solution ne sera pas unique, puisque en ajoutant une constante nous obtenons d'autres solutions. Cependant, deux solutions \`a cette \'equation  diff\'erent seulement d'une constante.  La  s\'erie de Fourier de $\phi$ 
\begin{equation*}
\mathcal{F}(\phi)(\theta)=\sum_{n\in \Z\setminus \{0\}}\hat{\phi}(n)e^{2\pi in\alpha}\quad, \quad \hat{\phi}(n)=\int_{\T}\phi(\theta)e^{-2\pi in\theta}d\theta
\end{equation*}
co\"{\i}ncide avec $\phi$ (on a convergence uniforme de toutes les d\'eriv\'ees). Soit $\psi$ une solution de classe $\cinf$ \`a l'\'equation (\ref{appppeqco}). On a donc
\begin{eqnarray*}
\mathcal{F}\big(\psi\circ R_{\alpha}-\psi\big)&=&\mathcal{F}(\phi)\\
\hat{\psi}(n)\big(e^{2\pi i n\alpha}-1\big)&=&\hat{\phi}(n)
\end{eqnarray*}
pour tout $n\in \Z\setminus \{0\}$.  Ceci nous sugg\`{e}re (et oblige) \`a d\'efinir une solution par la formule
\begin{equation}\label{appppserie1}
\psi(\theta)=\sum_{n\in \Z\setminus\{0\}}\frac{\hat{\phi}(n)}{e^{2\pi in\alpha}-1}e^{2\pi i n\theta}.
\end{equation}
Notons que la valeur de la moyenne  $\int_{\T}\psi(\theta)d\theta$ n'est pas fix\'ee \`a priori et nous choisissons la normalisation $\hat{\psi}(0)=0$. M\^eme si $\phi$ est de classe $\cinf$, et donc ses coefficients de Fourier d\'{e}croissent convenablement, la pr\'esence du facteur $\big(e^{2\pi in\alpha}-1\big)^{-1}$, qui peut devenir tr\`{e}s grand si $n\alpha$ s'approche d'un entier, peut produire que cette s\'erie ne corresponde pas \`a celle d'une fonction de classe $\cinf$. Plus pr\'{e}cis\'{e}ment on a l'in\'egalit\'e 
\begin{equation*}
4\|n\alpha\|\leq \big|e^{2\pi i n\alpha}-1\big|\leq 2\pi \|n\alpha\|.
\end{equation*}
On voit ainsi que les conditions arithm\'etiques de $\alpha$ apparaissent dans la discussion.
\begin{prop}
Si $\alpha\in \CD$ et $\phi$ est de classe $\cinf$ alors la s\'erie (\ref{appppserie1}) d\'{e}finit une solution de classe $\cinf$ \`a l'\'equation (\ref{appppeqco}). Si $\alpha\notin \CD$ alors il existe $\phi$ de classe $\cinf$ telle que la s\'erie (\ref{appppserie1}) n'est m\^eme pas une distribution $\quad_\blacksquare$
\end{prop}
 Soit $\beta\in \R$  et soit  $\phi:\T\to \C$ une fonction de classe $\cinf$. Dans ce travail nous devront consid\'erer  l'\'equation \emph{cohomologique tordue elliptiquement}
\begin{equation}\label{appppeqcot}
\psi(\tal)-e^{2\pi i\beta}\psi(\theta)=\phi(\theta).
\end{equation}
En appliquant la m\'ethode des s\'eries de Fourier nous obtenons la s\'erie
\begin{equation}\label{appppserie2}
\psi(\theta)=\sum_{n\in \Z}\frac{\hat{\phi}(n)}{e^{2\pi in\alpha}-e^{2\pi i\beta}}e^{2\pi i n\theta}.
\end{equation}
Nous voyons que cette fois les petits diviseurs qui apparaissent sont de la forme $\|n\alpha-\beta\|$, et ainsi les conditions arithm\'etiques sur la paire $(\alpha,\beta)$ entrent dans la discussion. Notons cependant que si la paire $(\alpha,\beta)$ est rationnellement independante, la s\'erie (\ref{appppserie2}) d\'efinit tous les coefficients de Fourier de la solution (ce qui ne se passait pas dans le cas de l'\'equation cohomologique classique pour le coefficient d'ordre $0$).
\begin{prop}\label{preliminares.contraejemplo}
Si $(\alpha,\beta) \in \CD_1$ et $\phi$ est de classe $\cinf$ alors la s\'erie (\ref{appppserie2}) d\'{e}finit une solution de classe $\cinf$ \`a l'\'equation (\ref{appppeqcot}). Si $\beta\notin \CD_1^{\alpha}$ alors il existe $\phi$ de classe $\cinf$ telle que la s\'erie (\ref{appppserie2}) n'est m\^eme pas une distribution. De plus, \'etant donn\'es $\eps>0, r\in \N$ on peut choisir $\phi$ de fa\c con  que sa taille $C^r$  soit plus petite que $\eps\quad_{\blacksquare}$
\end{prop}
  \section{D\'efinitions et pose du  probl\`eme}
  Dans ce travail nous allons supposer que les \emph{dhf} ainsi comme les courbes invariantes sont toujours de classe $\cinf$. Nous dissons qu'une courbe invariante $u:\T\to \C$ est \emph{indiff\'erente} si 
  \begin{equation*}
  \int_{\T}\log \big|\partial_{z}f\big(\theta,u(\theta)\big)\big|d\theta=0.
  \end{equation*} 
  Notons que les fonctions $f(\theta,\cdot)$ sont injectives et $\partial_zf(\theta,\cdot)$ ne s'annule donc pas. Nous traiterons das ce travail la persistance des courbes indiff\'erentes de degr\'e nul, c'est \`a dire, les courbes invariantes indif\'erentes pour lesquelles il se v\'erifie l'hypoth\`ese suivante: le degr\'e topologique de l'application
  \begin{equation*}
  \theta\longmapsto \partial_{z}f\big(\theta,u(\theta)\big)
  \end{equation*}
  est nul. De fa\c con \'equivalente, l'application ci-dessus est homotope \`a une constante dans $\C\setminus\{0\}$. Notons que dans ce cas l'application $\log \partial_{z}f\big(\theta,u(\theta)\big)$ est bien d\'efinie, $\mod 2\pi$.
  \paragraph{Le nombre de rotation transversal. } \`A une courbe invariante indiff\'erente de degr\'e nul nous lui associons un nombre qui mesure la vitesse moyenne d'enroulement des orbites proches autour de la courbe invariante:
      \begin{equation}\label{deffrn}
  \varrho_{tr}(u)=\frac{1}{2\pi i}\int_{\T}\log \partial_zf_{\theta}\big(u(\theta)\big)d\theta.
  \end{equation}
Cette quantit\'e est bien d\'efinie $\mod 1$ et nous  l'appelons le \emph{nombre de rotation transversal}.  Par exemple la \emph{dhf} donn\'e par  $F(\theta,z)=(\theta+\alpha,e^{2\pi i \beta}z)$ avec $\beta$ un nombre r\'eel dans l'interval $[0,1)$, est la plus simple mais pas la moins inter\'essante. La courbe $u=\{z\equiv 0\}_{\theta \in \T}$ est  lisse, invariante, de degr\'e nul,  indiff\'erente et $\varrho_{tr}(u)=\beta$. 
\paragraph{Forme Normale. } Soit $F$ une dynamique holomorphe fibr\'ee avec  une courbe invariante $\{u_0(\theta)\}_{\theta\in \T}$ indiff\'erente de classe $\cinf$,  de degr\'e nul et $\varrho_{tr}(u_0)=\beta\in \R$. Supposons que nous pouvons   r\'essoudre l'\'equation cohomologique
   \begin{equation}\label{introcohomologica2}
  \frac{u_1(\theta+\alpha)}{u_1(\theta)}e^{2\pi i \beta}=\partial_zf_{\theta}\big(u_0(\theta)\big)
  \end{equation}
avec $u_1$ une fonction \`a valeurs dans $\C$, qui ne s'annule pas, de degr\'e nul et de  classe $\cinf$. En faisant le changement de coordonn\'ees 
  \begin{equation*}\label{introH}
  H(\theta,z)=\big(\theta,u_0(\theta)+u_1(\theta)z\big)
  \end{equation*}
   nous obtiendrons une  forme normale pour $F$
   \begin{equation*}
(\theta,z)\longmapsto\big(\theta+\alpha,e^{2\pi i\beta}z+\rho(\theta,z)\big), 
\end{equation*}
o\`u la fonction $\rho(\theta,z)$ est d\'efinie dans le produit du cercle avec une voisinage de l'origin complexe, est de classe $\cinf$ et s'annule jusqu'au l'ordre $2$ en $z=0$.    On verra que l'\'equation (\ref{introcohomologica2}) peuve \^etre r\'esolue sous une hypoth\`ese  arithm\'etique diophantienne sur le nombre $\alpha$ (cf. section \ref{sec2}). 
    \newline
    
  La discussion pr\'ec\'edente nous permet donc de dire que la dynamique holomorphe fibr\'ee autour d'une courbe invariante peut \^etre vue, \`a la r\'esolution d'une \'equation cohomologique pr\`es,   comme une dynamique fibr\'ee par des germes holomorphes qui fixent l'origine, avec un nobre de rotation  bien pr\'ecis, au dessus d'une rotation irrationnelle du cercle. 
\paragraph{Petites perturbations. } Ce travail est consacr\'e \`a l'\'etude de la persistance d'une courbe invariante $u$ (lisse, de degr\'e nul et indiff\'erente), sous des petites perturbations sur la dynamique. D'apr\'es la discusion pr\'ec\'edente, sous l'hypoth\`ese diophantienne sur $\alpha$, nous dissons que $\tilde{F}$ est une petite perturbation de $F$ \emph{autour de la courbe} si, dans les coordonn\'ees de $H$ la transformation $\tilde{F}$ s'\'ecrit
\begin{equation*}
 (\theta,z)\longmapsto \big(\tal, \tilde{\rho}_0(\theta)+\big(\tilde{\rho_1}(\theta)+e^{2\pi i \varrho_{tr}({u})}\big)z+\tilde{\rho}(\theta,z)\big)
 \end{equation*}
 avec $\tilde{\rho}_0, \tilde{\rho}_1\, :\T\to \C$ des fonctions de classe $\cinf$ de taille petite (dans une topologie ad\'{e}quate), et avec $\tilde{\rho}(\theta,\cdot)$ holomorphe dans un disque $D_{r}$, $r>0$.  La taille de $\tilde{\rho}$ est comparable  \`a la taille de $\rho$ (dans un sens \`a pr\'eciser). Ces hypoth\`eses de proximit\'e entre la perturbation et la transformation originale seront bien pr\'ecis\'ees dans l'\'ennonc\'e du th\'eor\`me principal.  Notons que nous consid\'erons des perturbations qui changent seulement la partie holomorphe des dynamiques. Le nombre de rotation sur la base $\alpha$ est donc fix\'e.
 \paragraph{Familles \`a $1$ param\`etre. } Dans les r\'esultats de type KAM, une perturbation sur une dynamique elliptique donne lieu aussi \`a une perturbation sur les frequences concern\'ees, qui sont les valeurs qui contr\^olent la dynamique. Donc, nous ne pouvons pas esp\'erer d'obtenir les m\^emes ph\'enom\`enes dynamiques que dans la situation non perturb\'ee. Pour obtenir une persistance de ces ph\'enom\`enes  on introduit une famille \`a $1$ param\`etre de perturbations qui cherche \`a corriger les frequences \`a la valeur originale. De cette fa\c con on montre que, quite \`a corriger les frequences,  plusieurs propri\'et\'es dynamiques sont persistantes (lin\'earisation, existence d'objets invariants, etc.). On dit que la persistance est en co-dimension $1$ (voir \cite{BOST86}).  Dans notre travail nous perturbons seulement la partie holomorphe de la transformation, donc la frequence corresponde au nombre de rotation transversal, une donn\'ee de dimension $1$ complexe.
  \\  
  Soit  $F$ une \emph{fhd} et $u$ une courbe invariante indiff\'erente, avec $\varrho_{tr}(u)=\beta\in \R$. Soit  $\Sigma$ un ouvert de  $\C$.  Une  {petite perturbation transverse} ed $F$ est  une famille \`a $1$ param\`etre complexe  $\{F_t\}_{t\in \Sigma}$ de \emph{fhd} (une courbe complexe dans l'espace
   des \emph{fhd})  qui v\'erifie les hypoth\`eses suivantes:  chaque \'el\'ement  $F_t$ est une petite perturbation de  $F$ et le nombre de rotation transversal (m\^eme quand la courbe n'existe pas) bouge avec la famille $\{F_t\}_{t\in \Sigma}$ (voir le th\'eor\`eme  \ref{casoliso1} pour une d\'efinition pr\'ecise). Nous dissons que la courbe  $u$ est persistante si pour toute petite perturbation transverse  $\{F_t\}_{t\in \Sigma}$ il existe un param\`etre  $t^*\in \Sigma$ tel que  $F_{t^*}$ poss\`ede une courbe invariante indiff\'erente  $u^*$ avec $\varrho_{tr}(u^*)=\beta$. En gros, le principal r\'esultat de ce travail dit que, sauf pour une petite corrrection complexe, les courbes invariantes indiff\'erentes sont persistantes en classe $\cinf$ pourvu que les nombres de rotation v\'erifient une condition arithm\'etique de type diophantienne.  
  \\
\section{\'{E}nonc\'e du Th\'eor\`eme}\label{cinfinito}
 Nous consid\'erons une famille $\{F_s\}_{s\in \Sigma\subset \C}$ \`a un param\`etre $s\in \Sigma\subset \C$     complexe de \emph{dhf} o\`u chaque $F_s$ est une \emph{dhf} de classe $\cinf$.  Pour $\beta\in \R$ nous fixons la notation  $\lambda=e^{2\pi i \beta}$. Nous disons qu'une telle famille est \emph{lisse} si, \'ecrite sous la notation habituelle
\begin{eqnarray*}\label{2cinfnotacion}
F_s(\theta,z)&=&\big(\tal,f_s(\theta,z)\big)\\
&=&\big(\theta+\alpha,\rho_{0,s}(\theta)+\rho_{1,s}(\theta)z+\lambda z+ \rho_s(\theta,z)\big),
\end{eqnarray*}
les fonctions $(s,\theta,z)\mapsto \big(\rho_{0,s}(\theta), \rho_{1,s}(\theta), \rho_s(\theta,z)\big)$ sont des fonctions de classe $\cinf$. Les fibres $\rho_s(\theta,\cdot)$ sont holomorphes, continues jusqu'au bord $\partial \D$  et s'annulent  jusqu'au l'ordre 2 en $z=0$ pour pour tout $\theta$ dans $\T$ et tout $s$ dans $\Sigma$. Dans ce travail nous allons d\'esigner une famille de \emph{dhf} soit par $\{F_s\}$, soit par $\{f_s\}$ o\`u $f_s$ repr\'{e}sente la partie holomorphe de la dynamique $F_s$. Nous disons que $F_s$ est la \emph{dhf} associ\'ee \`a $f_s$. Nous donnons aussi l'adjectif de \emph{lisse} ($\cinf$) \`a une famille de parties holomorphes $\{f_s\}_{s\in\Sigma\subset \C}$. Avant d'\'{e}noncer le r\'esultat principal de ce travail nous allons introduire quelques notations.
\paragraph{Des notations et normes consid\'er\'ees. }  Nous allons \'ecrire $f'(\theta,z)$ au lieu de $\partial_z f(\theta,z)$ et $f^{(i)}(\theta,z)$ pour les d\'eriv\'ees d'ordre sup\'erieur $\partial_z^{i}f(\theta,z)$.
Pour une fonction $g:\T\to (B,|\cdot|_B)$ de classe $C^{\infty}$, o\`u $B$ est un espace de Banach avec norme $|\cdot |_B$, nous consid\'erons la norme $C^0$ et $C^r$ pour $r\in \N$ comme \'etant 
\[
\|g\|_0=\sup_{\theta\in \T}|g(\theta)|_B\quad,\quad \|g\|_r=\sum_{i=0}^r\Big \|\frac{\partial^ig}{\partial\theta^i}\Big \|_0.
\] 
 L'espace de Banach que nous aurons toujours en t\^ete sera l'espace $\OM$ des fonctions holomorphes du disque unit\'e complexe qui sont continues  jusqu'au bord. Pour une matrice  $A:\R^2\to \R^2$  nous allons noter
 \begin{equation*}
 \|A\|_{\LL}=\sup_{|v|=1}|Av|\quad,\quad  \big[A\big]_{\LL}=\inf_{|v|=1}|Av|
 \end{equation*} 
 sa norme usuelle et le \emph{plus petit valeur propre} respectivement,  o\`u $|\cdot|$ est une norme dans $\R^2$.  
 \begin{lema}\label{cinf2lemaA}
 Soit  $A:\R^2\to \R^2$ une matrice avec  $\big[A\big]_{\LL}\neq 0$. Alors la matrice $A$ est inversible et la norme de la matrice inverse $A^{-1}$ v\'erifie
 \begin{equation*}
 \|A^{-1}\|_{\LL}\leq \big[A\big]^{-1}_{\LL}  \quad_{\blacksquare}
 \end{equation*}
 \end{lema} 
 \begin{theo}\label{casoliso1}
    Pour toute paire $(\alpha,\beta)$ qui v\'erifie l'hypoth\`ese $\C\D_1$ et pour toutes constantes $L>1, M>1,T>1$ il existe $\tilde{\eps}>0$ qui d\'epend de $L, M, T, (\alpha,\beta)$,  un nombre entier naturel $r\geq 2$ qui d\'epend de la paire $(\alpha,\beta)$ et une constante positive universelle $C$ tels que, si une famille \`a un param\`etre complexe $\{f_s\}_{s\in \Sigma}$ de fonctions  de $\T$ vers $\OM$  v\'erifie   pour un certain $\eps$ dans $(0,\tilde{\eps}]$
    \begin{itemize}
    \item[$\bullet$] $\|\rho_{0,s}\|_r\leq \eps$ pour tout $s$ dans $D(0,2CL\eps)\subset\Sigma$
    \item[$\bullet$] $\|\rho_{1,s}\|_r\leq \eps$  pour tout $s$ dans $D(0,2CL\eps)\subset\Sigma$
    \item[$\bullet$] $\bigg[ \partial_s\int_{\T}\rho_{1,t}(\theta)d\theta\Big|_{s=0}\bigg]_{\LL}>L^{-1}$
    \item[$\bullet$] $\|\rho_s\|_r\leq M$ pour tout $s$ dans $D(0,2CL\eps)\subset\Sigma$
    \item[$\bullet$] $\|\partial_s^2\rho_{1,s}\|_0+\|\partial_s\rho_s\|_0\leq T$ pour tout $s$ dans $D(0,2CL\eps)$
    \end{itemize}
    alors il existe un param\`etre $s^*$ dans le disque $D(0,2CL\eps)$ et une courbe $u:\T\to \D$ de classe $\cinf$, indiff\'erente, de degr\'e nul, qui est invariante par la dynamique holomorphe fibr\'ee $F_{s^*}(\theta,z)=(\theta+\alpha,f_{s^*}(\theta,z))$,   et son nombre de rotation transversal est $\varrho_{tr}(u)=\beta$.
\end{theo}
Le th\'eor\`eme ci dessus d\'ecoulera d'une version plus faible qui cependant est plus adapt\'{e}e \`a la m\'ethode de d\'emonstration utilis\'ee
\begin{theo}\label{casoliso2}
Pour toute paire $(\alpha,\beta)$ qui v\'erifie l'hypoth\`ese diophantienne $\C\D_1$ il existe $\bar{\eps}>0$ et un nombre naturel $r\geq 2$ tels que si une fonction $f:\T\to \OM$ satisfait $\|f-\lambda z\|_r<\bar{\eps}$ alors il existe $t$ dans $\C$ et une courbe $u:\T\to \D$ de classe $\cinf$, indiff\'erente,  de degr\'e nul, qui est invariante par la dynamique holomorphe fibr\'ee $F^*(\theta,z)=(\theta+\alpha,e^tf)$,  et son nombre de rotation transversal est $\varrho_{tr}(u)=\beta$. 
\end{theo}
Le Th\'eor\`eme \ref{casoliso2} est donc la version du Th\'eor\`eme \ref{casoliso1} pour le cas particulier de la famille $\{f_t\}_{t\in \C}=\{e^tf\}_{t\in \C}$ et $f$ proche de $\lambda z$. 
Nous pouvons noter que la condition diophantienne demand\'ee dans les hypoth\`eses du th\'eor\`eme est exactement celle qui appara\^{\i}t comme n\'ecessaire et suffisante pour pouvoir toujours r\'esoudre l'\'equation lin\'earis\'ee associ\'ee au probl\`eme (voir section \ref{APPPcohomologique}). Cette simple observation nous permet de montrer que la condition diophantienne est optimale pour le probl\`eme de la persistance de la courbe invariante dans le cas $C^{\infty}$ comme le montre la 
\begin{prop}\label{contraejemplo}
Soit $\alpha\in \C\D$  et soit $\beta$ tel que la paire $(\alpha,\beta)$  ne satisfait pas la  condition diophantienne $\C\D_1$. Pour tous $\eps>0, r$ dans $\N$ il existe une fonction $a:\T\to \C$ de classe $C^{\infty}$ et de norme $C^r$ plus petite que $\eps$ telle que la famille \`a un param\`etre complexe $$F_t(\theta,z)=\big(\theta+\alpha,ta(\theta)+t\lambda z\big)$$ v\'erifie que pour tout $t$ dans $\C$ la dynamique holomorphe fibr\'ee $F_t$ ne poss\`ede  aucune courbe invariante de classe $C^{\infty}$ avec $\beta $ comme nombre de rotation transversal.
\end{prop}   
$\mathbf{Preuve.}$ Notons d'abord que pour une telle famille la d\'{e}riv\'{e}e par rapport \`a $z$ de la partie holomorphe est toujours \'egal \`a $t\lambda=te^{2\pi i \beta}$, donc le seul param\`etre qui permet l'existence  d'une courbe invariante avec nombre de rotation transversal \'egal \`a $\beta$ est $t=1$. L'hypoth\`ese de transversalit\'e est imm\'ediate. Supposons que pour une fonction $a:\T\to \C$ de classe $C^{\infty}$  nous avons une courbe invariante $u:\T\to \D$ pour la dynamique holomorphe fibr\'ee $F(\theta,z)=(\theta+\alpha,a(\theta)+\lambda z)$. Nous pouvons alors d\'eterminer de  fa\c con unique cette courbe \`a l'aide de l'\'equation de la courbe invariante
\begin{equation*}
 a(\theta)+e^{2\pi i\beta}u(\theta)=u(\theta+\alpha)
 \end{equation*}
 qui est bien une \'equation cohomologique comme celles \'etudi\'ees dans la section   \ref{APPPcohomologique}. La Proposition \ref{preliminares.contraejemplo} nous permet construire une fonction $a$ de fa\c con  que l'unique  courbe invariante solution \`a l'\'equation ci-dessus ne soit m\^eme pas une distribution $\quad_{\blacksquare}$
 \section{Le Th\'eor\`eme des Fonctions Implicites de Hamilton}
 Depuis les travaux de Sergeraert \cite{SERG72} et plus particuli\`erement  ceux de Herman \cite{HERM85},\cite{BOST86}, l'utilisation des th\'eor\`emes de fonctions implicites dans les espaces de Fr\'echet pour r\'esoudre des probl\`emes dynamiques faisant intervenir  des  petits diviseurs est devenue tr\`{e}s fructueuse. Cette technique repose sur le fait que la r\'esolution du  probl\`eme lin\'eaire associ\'e est fortement  reli\'e aux propri\'et\'es arithm\'etiques des fr\'{e}quences impliqu\'ees, comme nous pouvons  voir dans la discussion des \'equations cohomologiques (voir section  \ref{APPPcohomologique}). Un th\'eor\`eme de fonctions implicites assure en g\'{e}n\'{e}ral l'existence de solutions au le probl\`eme non lin\'eaire qui nous occupe pourvu que le probl\`eme lin\'eaire admette des solutions. Le Th\'eor\`eme des Fonctions Implicites de Hamilton nous permet d'appliquer cette technique dans le cadre des fonctions de classe $C^{\infty}$, qui d'habitude nous m\`{e}nent \`a travailler avec des espaces qui ne sont pas des espaces de Banach, mais des espaces de Fr\'{e}chet. La  diff\'erence essentielle avec les th\'eor\`emes des fonctions implicites classiques dans les espaces de Banach repose sur le fait que celui de Hamilton nous exige  de r\'esoudre le probl\`eme lin\'eaire associ\'e non pas seulement au  point o\`u la solution est connue, mais dans tout un voisinage de ce point. Les prochains paragraphes vont pr\'eciser tous les objets math\'ematiques qui interviennent dans l'\'{e}nonc\'{e} du Th\'eor\`eme de Hamilton ainsi que les espaces qui vont nous permettre de nous servir de ce th\'eor\`eme pour montrer la version faible du Th\'eor\`eme de la persistance de la courbe invariante dans le cas $C^{\infty}$ (Th\'eor\`eme \ref{casoliso2}). Nous renvoyons le lecteur aux articles \cite{HAMI82},\cite{BOST86} pour un traitement plus d\'etaill\'e du Th\'eor\`eme de Hamilton.
  \subsection{Bons Espaces de Fr\'echet}
 Nous disons que l'espace vectoriel
    topologique $E$ est un \emph{bon espace de Fr\'echet au sens de Hamilton} s'il existe une
    famille croissante de seminormes $\{\|\cdot \|_i\}_{i\in\N}$ qui d\'efinissent sa topologie,
    une famille d'op\'erateurs d'approximation et lissage
    ${\big\{S_t\big\}}_{t>1}$ et des constantes positives
    $C_{n,k}$ pour chaque paire $(n,k)$ dans $\N^2$ qui v\'erifient
    \begin{enumerate}
    \item $S_t:E\to E$ est une application lin\'eaire continue.
    \item Si $k\leq n ,\forall\ x\in E,\ \forall\ t\in\
    ]1,+\infty]\qquad \left\{ \begin{array}{l}
                        \|S_t(x)\|_n\leq C_{n,k}t^{n-k}\|x\|_k\\
                        \|[Id-S_t](x)\|_k\leq
                        C_{k,n}t^{k-n}\|x\|_n.
                        \end{array} \right. $
    \end{enumerate}
    Les derni\`{e}res in\'{e}galit\'{e}s impliquent des in\'egalit\'es de
    convexit\'e sur les seminormes $\|\ \|_i$, (\emph{Hadamard}): pour chaque paire $(n,k)$ dans $\N^2$ il existe des constantes positives $\tilde{C}_{n,k}$ tels que 
    \begin{equation}\label{hada}
     \|x\|_l\leq \widetilde{C}_{k,n}\|x\|_{k}^{1-\alpha}
    \|x\|_{n}^{\alpha}
    \end{equation}
     pour tout $x$ qui appartient \`a $E$, pour tous les entiers $k\leq l \leq n$,  o\`u $\alpha$ est d\'efini par $l=(1-\alpha)k+\alpha n$.

La somme directe (le produit)  $E\oplus F$ de deux bons espaces de Fr\'echet  est un bon espace de Fr\'echet avec les seminormes $\|(u,v)\|_{i_{E\oplus F}}=\|u\|_{i_E}+\|v\|_{i_F}$ et les op\'erateurs de lissage et approximation $S_t(u,v)=(S_t^Eu,S_t^Fv)$.
Nous disons que l'application
     $f:U\subset E\to F$ d'un ouvert $U$ d'un
     bon espace de Fr\'{e}chet $E$ vers un autre bon espace de Fr\'{e}chet
     $F$, est une \emph{bonne application au sens de Hamilton} si pour tout $x_0$ qui
     appartient \`a $U$ il existe un voisinage $V$ de $x_0$ dans
     $U$, un entier positif $r$  et pour tout $i\in \N$ des constantes positives $C_i$
      tels que
     \begin{equation*} \|f(x)\|_i\leq
     C_i(1+\|x\|_{i+r})
     \end{equation*}
     pour tout $x$ dans $V$ et pour tout $i$ dans $\N$. 
     \paragraph{Diff\'erentiabilit\'e au sens de G\^ateaux. } Soient $E$ et $F$ deux espaces  vectoriels topologiques, $U$ un ouvert de $E$ et $f$ une application de $U$ vers $F$. On dit que $f$ est de classe $C^1$ (au sens de G\^ateaux) lorsque
     \begin{itemize}
     \item[1. ] $f$ est continue.
     \item[2.  ] Il existe une application $Df:U\times E\to G$ continue, lin\'{e}aire en la deuxi\`eme coordonn\'ee et telle que pour tout $x$ dans $U$, $y$ dans $E$ on a 
     \[
     \lim_{t\to 0}\frac{1}{t}\Big\{f(x+ty)-f(x)\Big\}=Df(x)y.
     \]
     \end{itemize}
     Les applications de classe $C^k$ (au sens de G\^ateaux) sont d\'efinies par r\'ecurrence sur $k$: soit $k$ dans $\N\setminus\{0,1\}$, $f$ est dite de classe $C^k$ lorsqu'elle est de classe $C^1$ et que $Df$ est de classe $C^{k-1}$ sur l'ouvert $U\times E$ de $E\times E$.     Nous disons que $f$ est une bonne application de classe $C^k$ ($k$ dans $\N\cup\{\infty\}$) si $f$ est de classe $C^k$ au sens de G\^ateaux, et que$f$ ainsi que ses  d\'{e}riv\'{e}s jusque \`a l'ordre $k$ sont des bonnes applications (une telle d\'eriv\'e $D^if$ est une application \`a valeurs dans $F$ d\'efinie sur l'ouvert $U\times E^i$ de $E^{i+1}$).
    \begin{prop}\label{compo}
    Donnons nous trois bons espaces de Fr\'{e}chet $E,F,G$, $U$ un
    ouvert de $E$ et $V$ un ouvert de $F$. Si $f:U\mapsto V$ et
    $g:V\mapsto G$ sont des bonnes applications de classe $C^k$ ($k$ dans  $\N \cup \{\infty
    \}$), alors la composition $g\circ f$ est une bonne application de classe $C^k$. La projection $E\times F\longrightarrow E$ est une bonne application de classe $C^k$ pour tout $k$. 
    \end{prop}
            \subsection{Th\'eor\`eme des Fonctions Implicites}
    \begin{theo}[Fonction Implicite]\label{Hamil}
    Donnons nous trois bons espaces de Fr\'echet $E$, $F$, $G$, $U$ un
    ouvert de $E\times F$, $f:U\to G$ une bonne
    application de classe $C^r$ ($2\leq r\leq \infty$) et $(x_0,y_0)$ qui
    appartient \`a $U$. Supposons qu'il existe une voisinage $V_0$
    de $(x_0,y_0)$ et une bonne application continue et lin\'eaire
    dans la deuxi\`eme coordonn\'ee $L:V_0\times G\to F$
    telle que si $(x,y)$ appartient \`a $V_0$ alors $D_2 f(x,y)$ est
    inversible avec $L(x,y)$ comme son inverse.\
     On en d\'{e}duit alors que $x_0$ a un voisinage $W$ dans lequel est d\'efinie une
     bonne application de classe $C^r$, $g:W \to F$ telle que:
     \begin{enumerate}
     \item $g(x_0)=y_0$
     \item Pour tout $x$ dans $W$ la paire  $(x,g(x))$ appartient \`a $U$ et se v\'erifie que $
     f(x,g(x))=f(x_0,y_0)$
     \end{enumerate}
     En plus, si $x$ appartient \`a $W$, $y$ est dans un
     petit voisinage au tour de $y_0$ et on a que
     $f(x,y)=f(x_0,y_0)$ alors $y=g(x)$.
    \end{theo}
\subsection{Les bons espaces de Fr\'echet $\Gamma^\infty(\T,B)$} \label{frechet}
    Nous d\'efinissons $\Gamma^\infty(\T,B)$ comme l'espace des
    fonctions $f:\T\to B$ de classe $\cinf$ (au sens de G\^ateaux)  \`a valeurs
    dans un espace de Banach $(B,|\cdot|_B)$ et  nous le munissons de la famille de seminormes $C^r$.   
    Avec ces seminormes l'espace $\Gamma^\infty(\T,B)$ devient
    un espace de Fr\'echet. Nous d\'efinirons  des op\'erateurs
    d'approximation et lissage qui le font devenir   un bon espace de Fr\'echet
    sur lequel  nous pourrons appliquer le th\'eor\`eme
    \ref{Hamil}.
    \subsubsection{Op\'erateurs de lissage et approximation sur $\Gamma^\infty(\T,B)$}
        Les op\'erateurs que nous allons d\'efinir dans cette section sont classiques et nous renvoyons le lecteur aux articles \cite{HORM76},\cite{HERM83-2} pour les d\'emonstrations. La seule chose \`a souligner dans notre cas est que l'espace d'arriv\'ee des fonctions est un espace de Banach, o\`u on peut d\'efinir, de la m\^eme fa\c con que dans le cas r\'eel ou complexe,  l'int\'egrale de Riemann d'une fonction (voir \cite{KATZ04}).      L'op\'eration $\ast$ de convolution est ainsi bien d\'{e}finie. Par la suite nous    identifierons les fonctions qui appartiennent \`a      $\Gamma^\infty(\T,B)$ aux fonctions dans
        $\cinf_{\Z}(\mathbb{R},B)$, les fonctions de classe
        $\cinf $ qui sont $\Z$-p\'{e}riodiques \`a valeurs dans $B$. 
          Soit $\eta$ dans $\cinf (\mathbb{R},B)$ v\'erifiant  $\textrm{supp}(\eta)
        \subset [-1,1]$, $\eta(-x)=\eta(x)$ et $\eta(x)=1$ si $|x|\leq \frac{1}{2}$.
        Soit $\phi(x)=\int_{\R}e^{-2\pi i \xi
        x}\eta(\xi)d\xi$; on pose, pour   $t\geq1$,  $\phi_t(x)=t\phi
        (tx)$.
        Pour $f$ dans $\cinf_{\Z}(\mathbb{R},B)$ on d\'efinit
        l'op\'erateur de lissage et approximation
        \begin{equation}\label{oplisage}
        S_tf=f\ast\phi _t =
        \int_{\R}f(x-y)\phi_t(y)dy\quad \in \quad
        \cinf_{\Z}(\mathbb{R},B).
        \end{equation}
        Par la formule d'inversion de Fourier, on a pour tout $v$ dans $B$ et pour tout $n$ dans $\Z$
        \begin{equation*}
        S_t(ve^{2\pi i n\theta})=v\eta\Big(-\frac{n}{t}\Big)e^{2\pi i n\theta}
        \end{equation*}
        donc $S_tf$ est un polyn\^ome trigonom\'{e}trique de degr\'e au plus $|t|$. Les op\'erateurs $S_t$ ont les
        propri\'{e}t\'{e}s suivantes de lissage et approximation 
        \begin{prop}
            Pour chaque paire $(k,n)$ dans $\N^2$ il existe  des
            constantes positives $C_{k,n}$ telles que si $f$
            appartient \`a $\Gamma^\infty(\T,B)$ et $n\geq k$ alors pour tout $t\geq1$
        \begin{enumerate}
        \item $\|S_t f\|_n\leq C_{k,n}t^{n-k}\|f\|_k$.
        \item $\|S_t f -f\|_k \leq C_{n,k}t^{k-n}\|f\|_n$.
        \end{enumerate}
        \end{prop}
        \section{Preuve du Th\'eor\`eme \ref{casoliso2}}
        \`A partir des bons espaces de Fr\'echet $\Gamma^{\infty}(\T,B)$
        d\'efinis dans la section \ref{frechet} on d\'efinit ici les
        bons espaces de Fr\'echet que nous allons utiliser lors de  la
        preuve du Th\'eor\`eme \ref{casoliso2} :  
        \begin{enumerate}
        \item L'espace $\cinf(\T,\C)$ des fonctions de classe
        $\cinf$ de $\T$ vers  l'espace de Banach $\C$ est le bon espace de Fr\'{e}chet
        $\Gamma^{\infty}(\T,\C)$.
                \item Le bon espace de Fr\'{e}chet  $\gaminf$ est l'espace des fonctions de classe
        $\cinf$ de $\T$ vers l'espace de Banach $\OM$. Ceci est l'espace des parties holomorphes des dynamiques holomorphes
        fibr\'ees de classe $\cinf$. Rappelons  que l'espace de Banach
        $\OM$ est l'espace des fonction holomorphes sur $\D$ qui
        sont continues  sur jusqu'au bord, avec la norme
        $|f|_{\OM}=\sup_{|z|<1}|f(z)|$.
       \end{enumerate}
       \begin{lema}\label{fronu2}
       Soit l'ensemble ouvert $A= \{u\in \Cinf\  \big |\  \|u\|_0<1/2\}$ de $\Cinf$. Les applications 
       \begin{itemize}
       \item[1. ] \begin{eqnarray*}
       \gaminf &\longrightarrow& \Gamma^{\infty}\big(\T,\Omega(\overline{D(0,3/4)})\big)\\
       f&\longmapsto& f^{(i)}=\frac{\partial^if}{\partial z^i},
       		\end{eqnarray*}
		\item[2. ]\begin{eqnarray*}
		\Cinf \times \Cinf &\longrightarrow& \Cinf\\
		(u,v)&\longmapsto&uv\\
		(u,v)&\longmapsto&u+v,
		\end{eqnarray*}
		\item[3. ]\begin{eqnarray*}
		\Cinf_0&\longrightarrow& \Cinf_0\\
		v&\longmapsto&v^{-1}\\
		v&\longmapsto&e^v
		\end{eqnarray*}
		o\`u $\Cinf_0$ est l'espace ouvert de $\Cinf$ des fonctions qui ne s'annulent pas, 
		\item[4. ]\begin{eqnarray*}
		\Cinf_{d_0}&\longrightarrow&\Cinf\\
		h&\longmapsto&\log h
		\end{eqnarray*}
		o\`u $\Cinf_{d_0}$ est l'espace ouvert de $\Cinf$ des fonctions de degr\'e nul, 
		\item[5. ]\begin{eqnarray*}
		\Cinf&\longrightarrow&\C\\
		v&\longmapsto& \int_{\T}v(\theta)d\theta ,
		\end{eqnarray*}
		\item[6. ] l'application $\eta_r$
		\begin{eqnarray*}
		\Gamma^{\infty}\big(\T,\Omega(\overline{D(0,r)})\big)\times A&\stackrel{\eta_r}{\longrightarrow}& \Cinf\\
		(f,u)&\longmapsto& f\big(\cdot, u(\cdot)\big)
		\end{eqnarray*}
		pour tout $r$ dans $(1/2,1]$,
       \end{itemize}
       sont de bonnes applications au sens de Hamilton.
       \end{lema}
       $\Pru $ Le point $1. $ est une cons\'{e}quence des estimations de Cauchy. Les points $2., 3., 4., 5., $ sont classiques et en fait ce sont des bonnes applications de classe $\cinf$. Nous montrerons par la suite le point $6$. Rappelons  que le fait d'\^etre une bonne
            application est un fait local, donc \'etant donn\'ee
            $(\overline{f},\overline{u})$ dans $\gaminf\times A$ nous  fixons un voisinage born\'e
            $C^0$ de $\overline{f}$ dans $\Gamma^{\infty}(\T,\Omega(\overline{D(0,r)}))$ et un voisinage $C^0$ de
            $\overline{u}$ dans $A\subset \Cinf $.            Notons que les estimations de Cauchy assurent
            l'existence pour chaque paire $(i,j)$ dans $\N^2$, de constantes positives $C_{i,j}$ telles
            que \begin{equation}\label{cauchy}
            \Big\|\frac{\partial^j}{\partial\theta^j}f^{(i)}(\theta,u(\theta))\Big
            \|_0\leq C_{i,j}\Big \|\frac{\partial^j f}{\partial \theta^j}\Big
            \|_0 \end{equation}
            car nous pouvons permettre des pertes de rayon uniformes.            Calculons quelques d\'eriv\'ees de $\eta_r$ par rapport
            \`a $\theta$:
           \begin{eqnarray*}
            \frac{\partial\eta_r(f,u)}{\partial\theta}&=&\frac{\partial
            f}{\partial\theta}(\cdot,u(\cdot))+f'(\cdot,u(\cdot))\frac{\partial u}{\partial
            \theta} \\
             \frac{\partial^2\eta_r(f,u)}{\partial\theta^2}&=& \frac{\partial^2f}{\partial\theta^2}(\cdot,u(\cdot))
            + \frac{\partial f'}{\partial\theta}(\cdot,u(\cdot))\frac{\partial u}{\partial\theta}
            + \frac{\partial f'}{\partial\theta}(\cdot,u(\cdot))\frac{\partial
            u}{\partial\theta}\\
            &  & + f''(\cdot,u(\cdot))\Big (\frac{\partial u}{\partial\theta}\Big )^2
            + f'(\cdot,u(\cdot))\frac{\partial^2
            u}{\partial\theta^2}.\\
            \end{eqnarray*}
           Plus g\'{e}n\'{e}ralement  nous voyons que gr\^{a}ce aux estimations de Cauchy (\ref{cauchy}) il nous
            suffit d'estimer les termes de la forme 
            \begin{equation*}
            \Big\| \frac{\partial^j
            f}{\partial\theta^j}\Big (\frac{\partial u}{\partial\theta}\Big )^{i_1}\Big (\frac{\partial^2 u}{\partial\theta^2}\Big )^{i_2}\dots
            \Big (\frac{\partial^n u}{\partial\theta^n}\Big )^{i_n}\Big\|_0
          \end{equation*} 
          avec $j+i_1+2i_2+\ldots +ni_n=n$, par rapport aux
          seminormes $\|f\|_n$ et $\|u\|_n$.
          En utilisant les in\'{e}galit\'{e}s de convexit\'e de Hadamard (\ref{hada})
          nous avons les estimations 
                \begin{eqnarray*}
                    \Big\| \frac{\partial^j
                    f}{\partial\theta^j}\Big\|_0&\leq&
                    \tilde{C}^{\Gamma^{\infty}}_{0,n}\|f\|_0^{\frac{n-j}{n}}\|f\|_n^{\frac{j}{n}}\leq A_n\|f\|_n^{\frac{j}{n}}\\
                    \Big \| \frac{\partial^s u}{\partial\theta^s}\Big
                    \|_0^{i_s}&\leq&\big
                    (\tilde{C}^{C^{\infty}}_{0,n}\|u\|_0^{\frac{n-s}{n}}\|u\|_n^{\frac{s}{n}}\big
                    )^{i_s}\leq B_n\|u\|_n^{\frac{si_s}{n}}
                \end{eqnarray*}
                pour des constantes positives $A_n,B_n$ qui d\'ependent seulement du voisinage de $\bar{f}$ \`a part de $n$.
               En  consid\'erant un produit ad\'{e}quat des in\'{e}galit\'{e}s ci-dessus
                nous avons
                \begin{equation*} \Big \|\frac{\partial^j
            f}{\partial\theta^j}(\theta,u(\theta))\Big (\frac{\partial u}{\partial\theta}\Big )^{i_1}\Big (\frac{\partial^2 u}{\partial\theta^2}\Big )^{i_2}\dots
            \Big (\frac{\partial^n u}{\partial\theta^n}\Big
            )^{i_n}\Big \|_0\leq
            D_n\|f\|_n^{\frac{j}{n}}\|u\|_n^{\frac{n-j}{n}}
            \end{equation*}
            pour des constantes positives $D_n$. 
            L'in\'egalit\'e de Young implique finalement que 
            \begin{equation*}\label{youn}
            \Big \|\frac{\partial^j
            f}{\partial\theta^j}(\theta,u(\theta))\Big (\frac{\partial u}{\partial\theta}\Big )^{i_1}\Big (\frac{\partial^2 u}{\partial\theta^2}\Big )^{i_2}\dots
            \Big (\frac{\partial^n u}{\partial\theta^n}\Big
            )^{i_n}\Big \|_0\leq D_n\big
            (\|f\|_n+\|u\|_n\big).
            \end{equation*}
            On peut en d\'eduire donc que l'application $\eta_r$ est une bonne application de classe $C^0\quad_{\blacksquare}$
            \begin{coro}\label{fronu}
                    L'application $\eta_r$ (cf. 6. du Lemme \ref{fronu2})  est une bonne application de classe $C^{\infty}$ pour tout $r$ dans $(1/2,1]$.
            \end{coro}   
            $\Pru $ Nous devons montrer que toutes les d\'{e}riv\'{e}es
            $D^i\eta_r$ sont des bonnes applications de classe $C^0$. 
            Voyons qu'est ce qui se passe avec $D\eta_r$
            \begin{eqnarray*}
             D\eta_r(f,u)(\Delta f,\Delta
            u)&=& \lim_{t\to 0}\frac{1}{t}\big [(f+t\Delta
            f)(\cdot,u+s\Delta u(\cdot))-f(\cdot,u(\cdot))\big ]\\
             &=&f'(\cdot,u(\cdot))\Delta u+\Delta
            f(\cdot,u(\cdot)).
            \end{eqnarray*}
         Ceci est une fonction \`a $4$ variables dans un bon espace de Fr\'echet (le produit de bons espaces) qui peut s'\'ecrire comme la composition de bonnes applications qui appara\^{\i}ssent dans le lemme pr\'{e}c\'{e}dent, donc c'est une bonne application. Pour les d\'eriv\'ees d'ordre sup\'{e}rieur il se passe de fa\c con analogue $\quad_\blacksquare$ 
                    \subsection{Les bons op\'{e}rateurs $\M_{\alpha}, \M_{\alpha,\beta}$}\label{cinf2emes}
            Soit $v$ dans $\Cinf$. Nous consid\'erons les \'equations cohomologiques
            \begin{eqnarray}\label{cinf2coho1} 
            \phi(\theta)-\phi(\tal)=v(\theta)\\\label{cinf2coho2}
            \lambda\tilde{\phi}(\theta)-\tilde{\phi}(\tal)=v(\theta).
            \end{eqnarray}
            Si $\int_{\T}v(\theta)d\theta=0$ nous avons d\'ej\`a vu que sous l'hypoth\`ese diophantienne pour $\alpha$ l'\'equation (\ref{cinf2coho1}) a pour solution 
            \begin{equation*}\label{cinf2solucion1}
            \phi(\theta)=\sum_{n\in \Z}\frac{\hat{v}(n)}{e^{2\pi i n\alpha}-1}e^{2\pi in\theta}
            \end{equation*}
            avec la normalisation $\hat{\phi}(0)=0$. Ainsi $\phi$ est de classe $\cinf$. De m\^eme, sous l'hypoth\`ese diophantienne pour la paire $(\alpha,\beta)$ la solution de l'\'equation (\ref{cinf2coho2}) est
             \begin{equation*}\label{cinf2solucion2}
            \tilde{\phi}(\theta)=\sum_{n\in \Z}\frac{\hat{v}(n)}{e^{2\pi i n\alpha}-e^{2\pi i \beta}}e^{2\pi in\theta}
            \end{equation*}
            qui est aussi de classe $\cinf$. Nous d\'efinissons les op\'{e}rateurs $\C-$lin\'eaires
            \begin{eqnarray*}
            \M_{\alpha}:\Cinf_*&\to&\Cinf\\
            v&\mapsto&\phi\\
            \M_{\alpha,\beta}:\Cinf&\to&\Cinf\\
            v&\mapsto&\tilde{\phi}
            \end{eqnarray*}
          o\`u $\phi$ et  $\tilde{\phi}$ sont d\'efinies ci-dessus et $\Cinf_*$ est l'espace des fonctions dans $\Cinf$ de moyenne nulle. Ainsi d\'efinis ces op\'{e}rateurs sont inversibles (compte tenue de la normalisation $\int_{\T}\phi=0$).
          \begin{lema}
          Les op\'{e}rateurs $\M_{\alpha}, \M_{\alpha,\beta}$ sont des bonnes applications au sens de Hamilton.
          \end{lema}
            $\Pru $ La th\'{e}orie classique des s\'eries de Fourier d'une fonction $f$ dans $\Cinf$ assure que les coefficients de sa s\'erie de Fourier v\'erifient les in\'egalit\'es
              \begin{eqnarray}\label{cinf2estrella1}
               \|f\|_i&\leq& B_i \sup_{k\in \mathbb{Z}}\Big
                                ((1+|k|)^{i+2}|\hat{f}(k)|\Big)\\\label{cinf2estrella2}
                                      \sup_{k\in
                                    \mathbb{Z}}\Big ((1+|k|)^{i}|\hat{f}(k)|\Big)&\leq&C_i\|f\|_i
                                                   \end{eqnarray}
pour des constantes positives $B_i,C_i$  qui d\'ependent seulement de $i$ dans $\N$. Soit $v$ dans $\Cinf$. Les in\'{e}galit\'{e}s ci-dessus et la condition diophantienne sur $\alpha$ impliquent que pour tout $i\geq 0$ il existe des constantes positives $C'_i, B'_{i+2}$ tels que
  \begin{eqnarray*} 
  \|\M_{\alpha}(v)\|_i&\leq& C_i \sup_{k\in \mathbb{Z}}\Big
                                ((1+|k|)^{i+2}\Big|\frac{\hat{v}(k)}{e^{2\pi in\alpha}-1}\Big|\Big)\\
                                    &\leq& C'_i \sup_{k\in
                                    \mathbb{Z}}\Big ((1+|k|)^{i+4}|\hat{v}(k)|\Big)\\
                                    &\leq& B'_{i+2}\|v\|_{i+2}.
                \end{eqnarray*}
La preuve pour $\M_{\alpha,\beta}$ \'etant analogue nous l'omettons $\quad_{\blacksquare}$   
            \subsection{Correction du nombre de rotation transversal}\label{correccion}
Pour toute $f$ dans $\gaminf$ telle que la taille  $C^0$ de $f-\lambda z$ est suffisamment petite et pour toute courbe continue $u:\T\to \D$  aussi de taille $C^0$ petite nous pouvons d\'efinir uniform\'{e}ment la fonction  $\log f'(\theta,u(\theta))$ qui v\'erifie $\log \lambda=2\pi i \beta$, car le degr\'e de la fonction $f'(\theta,u(\theta))$ d\'epend contin\^{u}ment de la paire $(f,u)$.  Dans ce cas nous pouvons calculer l'int\'egrale suivante 
        \begin{equation}\label{integ}
            \mathcal{I}(f,u)=\frac{1}{2\pi i}\int \log
            f'(\theta,u(\theta))d\theta .
                    \end{equation}  
                    Notons que quand la courbe $u$ est invariante par la dynamique associ\'ee \`a $f$, et de degr\'e nul, ce nombre co\"{\i}ncide avec le nombre de rotation transversal $\varrho_{tr}(u)$.
                    Avec cette d\'efinition nous pouvons donner du sens \`a l'affirmation suivante  ``pour n'importe  quelle paire $(f,u)$  il existe toujours un nombre complexe $t$ qui corrige le nombre de rotation transversal \`a la valeur $\beta$''.
            \begin{prop}\label{condiint} Pour tout nombre r\'eel $\beta$ dans $[0,1)$, pour toute application $f$ dans un voisinage $U$ de $f_0\equiv \lambda z$  dans $\gaminf$ et pour toute courbe $u$ dans une voisinage $V$ de $u_0\equiv 0$ dans $\Cinf$, la fonction $t:U\times V\to \C$ d\'efinie  par l'\'egalit\'e $\mathcal{I}(e^tf,u)=\beta$ est une bonne application de classe $\cinf$. En plus, il existe une constante positive universelle $C$ telle que  
            \begin{equation}\label{estimacionT}
            |t|\leq C\|f'-\lambda\|_{C^0}.
            \end{equation}
                           \end{prop}
            $\Pru $ En fait, la formule (\ref{integ}) nous permet calculer explicitement la valeur de $t$. Plus pr\'{e}cis\'{e}ment, soient $U,V$ les voisinages qui nous permettent de calculer l'int\'egrale $\mathcal{I}$ comme au paragraphe pr\'{e}c\'{e}dent, alors
            \begin{eqnarray}\nonumber
            \mathcal{I}(e^tf,u)&=&\frac{1}{2\pi i}\int \log
            e^tf'(\theta,u(\theta))d\theta \\\label{formulat}
            \beta&=&\frac{1}{2\pi i}\Big(t+\int \log
            f'(\theta,u(\theta))d\theta   \Big).         
                        \end{eqnarray}
                        Ainsi $t$ est une bonne application de classe $\cinf$ d'apr\'es la 
                         Proposition \ref{compo} et les Lemmes \ref{fronu2}, \ref{fronu}. L'estimation (\ref{estimacionT}) s'obtient facilement \`a l'aide de  (\ref{formulat})$\quad_\blacksquare$
            
            Dor\'enavant chaque fois que nous \'ecrivons $e^tf(\theta,u(\theta))$, et qu'il n'y a pas lieu  \`a confusion, il faudra penser toujours que $t=t(f,u)$.
                        \subsection{La bonne application $u_1$}\label{equationutil}          
                         L'\'egalit\'e  $\mathcal{I}(e^tf,u)=\beta$ et l'hypoth\`ese $\alpha$ dans $\C\D$ nous permettent  de
            r\'esoudre l'\'equation 
           \begin{equation}\label{cinf2coomouno} 
           e^tf'(\theta,u(\theta))=e^{2\pi
            i \beta}\frac{u_1(\theta+\alpha)}{u_1(\theta)}
            \end{equation}
             avec $u_1 \in \Cinf_0$, et o\`u en plus  $(f,u)\mapsto u_1$  est une bonne application au sens de Hamilton. En effet, dans les bons voisinages, o\`u la fonction $t$ est bien d\'efinie, l'\'equation
             (\ref{cinf2coomouno}) se r\'eduit  \`a l'\'equation cohomologique
             \begin{equation*} 
             t+\log f'(\theta,u(\theta))= 2\pi i\beta
             +\tilde{u}_1(\theta+\alpha) -\tilde{u}_1(\theta).
             \end{equation*}
             On voit que la condition $\mathcal{I}(e^tf,u)=\beta$ est la condition int\'{e}grale n\'ecessaire pour r\'esoudre cette \'equation. 
             Pour assouplir la notation nous introduisons la fonction $l(f,u)(\theta)=\log f'(\theta,u(\theta))$. Il est clair que l'application $(f,u)\mapsto l(f,u)$ est une bonne application dans les voisinages que nous consid\'erons ici. Nous avons alors
                       \begin{equation*}
                       \tilde{u_1}(\theta)=\M_{\alpha}(l(u,f))(\theta).
                                    \end{equation*}
                               L'application
                    $(f,u)\mapsto \tilde{u_1}$ est une bonne application
                    au sens de Hamilton d'apr\`es les r\'esultats de la section \ref{cinf2emes}, et il en est de m\^eme de l'application 
                    \begin{equation*}
                    (f,u)\longmapsto u_1=e^{\tilde{u}_1}.
                    \end{equation*}
                    Si en autre, si les voisinages $V,U$ sont suffisamment petits la distance $|u_1(\theta)-1|$ est tr\`{e}s petite et uniform\'{e}ment born\'ee, et en particulier $u_1(\theta)$ ne s'annule pas. Il est imm\'{e}diat de v\'erifier que $u_1$ satisfait l'\'equation (\ref{cinf2coomouno}).
\subsection{Le Th\'eor\`eme \ref{casoliso2} sous la  forme d'un probl\`eme de fonction implicite}
             Consid\'erons l'application  $\Theta :\Big(U\subset \gaminf\Big) \times \Big(V\subset \cinf(\T,\C)\Big) \to
             \cinf(\T,\C)$ d\'efinie par
             \begin{equation}\label{teta}\Theta(f,u)=
             e^tf(\theta,u(\theta))-u(\theta+\alpha)
             \end{equation}
             o\`u      l'application $t(f,u)$ et les ensembles $U$,$V$ sont d\'efinis
             dans  la Proposition \ref{condiint}. Cette application est une bonne application de classe $C^{\infty}$.       On peut sans peine voir que le fait que $u$ soit
                invariante pour la dynamique holomorphe fibr\'ee associ\'ee \`a  $e^tf$ est \'equivalente
                au fait que $\Theta(f,u)$ soit identiquement
                nulle.      Notons que dans le cas o\`u  $f_0\equiv \lambda z$ et $u_0\equiv 0$ nous avons 
             \begin{equation*}
             \Theta (f_0,u_0)=0.
             \end{equation*}
             La preuve du Th\'eor\`eme \ref{casoliso2}
             sera une application directe du Th\'eor\`eme des Fonctions
             Implicite de Hamilton, c'est-\`a-dire, la courbe
             invariante $u$ sera d\'efinie d'une fa\c con implicite
             \`a partir de l'\'equation $\Theta(f,u)=0$ autour de la solution d\'ej\`a connue $(f_0,u_0)$.
             \subsubsection{Inversion de la diff\'erentielle}
                La partie la plus importante pour appliquer le Th\'eor\`eme de Hamilton est l'inversion de la diff\'erentielle par rapport \`a la deuxi\`eme variable, ce qui dans notre cas se traduit par trouver une bonne application continue $L(f,u,\Delta g)$ , lin\'eaire en $\Delta g$, d\'efinie pour toute paire $(f,u)$ dans un voisinage de $(f_0,u_0)$ et pour tout $\Delta g$ dans $\Cinf$ de fa\c con  que si nous posons $\Delta u=L(f,u,\Delta g)$ nous avons
                                \begin{equation}\label{cinf2D2}
                                D_2\Theta(f,u)\Delta u = \Delta g.
                                 \end{equation}                              
                                 Dans ce qui se suit nous allons r\'esoudre cette \'equation (en $\Delta u$ dans $\Cinf$) d'une fa\c con formelle.  
             La diff\'erentielle partielle de $\Theta$ par rapport \`a
             $u$ en la direction $\Delta u$ est 
             \begin{equation*}
             D_2\Theta(f,u)\Delta u =e^tf(\cdot,u(\cdot))\big(\partial_ut\cdot\Delta u\big)+e^tf'(\cdot,u(\cdot))\Delta u-\Delta u(\cdot+\alpha)
                       \end{equation*}
o\`u $\partial_ut\cdot\Delta u$ est la diff\'erentielle de $t$ par rapport \`a $u$ dans la direction $\Delta u$, qui est un nombre complexe. Nous prenons $u_1(f,u)=u_1$ comme dans  (\ref{cinf2coomouno}). L'\'equation (\ref{cinf2D2}) devient donc 
             \begin{equation*}
              e^tf(\cdot,u)\big(\partial_ut\cdot\Delta u\big)+\frac{u_1(\cdot+\alpha)}{u_1}e^{2\pi i
             \beta}\Delta u
             -\Delta u(\cdot + \alpha)=\Delta g .
             \end{equation*}
             Avec la notation $\widetilde{\Delta u}=\frac{\Delta u}{u_1}$, 
              $\widetilde{\Delta g}=\frac{\Delta g}{u_1(\cdot+\alpha)}$, 
             $\tilde{f}=\frac{e^tf(\cdot,u)}{u_1(\cdot+\alpha)}$  
              nous avons 
             \begin{equation*}
             e^{2 \pi i \beta}\widetilde{\Delta u}(\theta)-\widetilde{\Delta
             u}(\theta+\alpha)=\widetilde{\Delta g}(\theta)-\big(\partial_ut\cdot\Delta
             u\big)\tilde{f}(\theta).
             \end{equation*}
             En appliquant l'op\'erateur $\C-$lin\'eaire (inversible) $\M_{\alpha,\beta}$ aux deux cot\'es on obtient de fa\c con \'equivalente
           \begin{equation}\label{cinf2otraulo}
           \widetilde{\Delta u}=\M_{\alpha,\beta}(\widetilde{\Delta g})-\big(\partial_ut\cdot\Delta u\big)\M_{\alpha,\beta}(\tilde{f}).
           \end{equation}
           On peut calculer explicitement la valeur de $\partial_ut\cdot\Delta u$, car en d\'erivant l'\'egalit\'e $\mathcal{I}(e^tf,u)=\beta$ dans la direction $\Delta u$ nous obtenons
           \begin{equation*}
           \frac{1}{2\pi i}\bigg(\big(\partial_ut\cdot\Delta u\big)+ \int_{\T}\frac{f''(\theta,u(\theta))\Delta u(\theta)}{f'(\theta,u(\theta))}d\theta\bigg)=0
           \end{equation*} 
           ce qui ajout\'e \`a l'\'egalit\'e (\ref{cinf2otraulo}) nous donne
           \begin{eqnarray}\nonumber
           \partial_ut\cdot\Delta u&=& -\int_{\T}\frac{f''(\theta,u(\theta))\Delta u(\theta)}{f'(\theta,u(\theta))}d\theta\\\nonumber
           &=&-\int_{\T}\frac{f''(\theta,u(\theta))u_1(\theta)\Big(\M_{\alpha,\beta}(\widetilde{\Delta g})(\theta)-\big(\partial_ut\Delta u\big)\M_{\alpha,\beta}(\tilde{f})(\theta)\Big)}{f'(\theta,u(\theta))}d\theta\\\label{425}
           \partial_ut\cdot\Delta u&=&\frac{-\int_{\T}\frac{f''(\theta,u(\theta))u_1(\theta)}{f'(\theta,u(\theta))}\M_{\alpha,\beta}(\widetilde{\Delta g})(\theta)d\theta}{1-\int_{\T}\frac{f''(\theta,u(\theta))u_1(\theta)}{f'(\theta,u(\theta))}\M_{\alpha,\beta}(\tilde{f})(\theta)d\theta}.
                      \end{eqnarray}
                      Nous posons finalement 
 \begin{equation*}
 \Delta u=u_1\Big(\M_{\alpha,\beta}(\widetilde{\Delta g})-\big(E\M_{\alpha\beta}(\tilde{f})\Big)
 \end{equation*}
 o\`u $E=\partial_u t\cdot \Delta u$ comme dans (\ref{425}).   
 On voit que si $f''(\theta,u(\theta))$ est suffisamment petit par rapport \`a  la taille $C^0$ de $u_1$ et $\M_{\alpha,\beta}(\tilde{f})$, le nombre complexe $\partial_ut\Delta u$ est bien d\'efini et en plus l'application $(f,u,\Delta g)\mapsto \partial_ut\Delta u$ est une bonne application au sens de Hamilton.  Il  est direct aussi que la d\'efinition ci-dessus pour  $\Delta u$ v\'erifie l'\'equation (\ref{cinf2D2}), et l'application $(u,f,\Delta g)\mapsto \Delta u$ est une bonne application au sens de Hamilton dans un voisinage ad\'equat de la paire $(f_0,u_0)$. Or, les conditions sur $f'', u_1, M_{\alpha,\beta}(\tilde{f})$  s'obtiennent en r\'etrecissant suffisamment le voisinage autour de $(f_0,u_0)$ qui sert \`a calculer $t$. Nous avons en main tous les ingr\'edients pour appliquer le Th\'eor\`eme de Hamilton \`a l'application $\Theta$: Il existe une bonne application de classe $\cinf$
 \begin{eqnarray*}
  \Xi:\widetilde{U}\subset \gaminf&\longrightarrow& \C\times\Cinf\\
  \Xi(f)&\longmapsto&(t,u)
  \end{eqnarray*}
   d\'efinie dans un voisinage $\tilde{U}$ de $f_0$ (une $C^r$ boule, pour un certain $r\in \N$)   telle que  la courbe $u$ est invariante pour la dynamique holomorphe fibr\'ee associ\'ee \`a $e^tf$, de degr\'e nul et de nombre transversal de rotation \'egal \`a $\beta$, ce qui donne une version plus pr\'ecise du Th\'eor\`eme \ref{casoliso2}. Soient $\eps>0$ et $r\in \N$. Nous posons  
   \[
   B_{\bar{\eps}}^r=\{f \in \gaminf \ \big| \ \|f-\lambda z\|_r<\bar{\eps}\}.
   \]
  \begin{theo}\label{casoliso3}
  Pour toute paire $(\alpha,\beta)$ qui v\'erifie l'hypoth\`ese $\C\D_1$ il existe $\bar{\eps}>0$, un nombre naturel $r\geq 2$ et une bonne application de classe $\cinf$ au sens de Hamilton, $\Xi:B_{\bar{\eps}}^r\subset \gaminf \to \C\times \Cinf$ tels que 
  \begin{itemize}
  \item[$\bullet$] $\Xi(\lambda z)=(0,0)$.
  \item[$\bullet$] Si on \'ecrit $\Xi(f)=(t,u)$ alors la courbe $u$ est invariante par la dynamique holomorphe fibr\'ee $F^*(\theta,z)=(\theta+\alpha,e^tf(\theta,z))$, est de degr\'e nul  et son nombre de rotation transversal est $\varrho_{tr}(e^tf,u)=\beta$.
  \item[$\bullet$] La projection sur la premi\`ere coordonn\'ee $\pi_1\Xi(f)$ co\"{\i}ncide avec  $t(f,\pi_2\Xi(f))$ o\`u $\pi_2\Xi(f)$ est la projection sur la deuxi\`{e}me coordonn\'ee et $t$ est la fonction d\'efinie dans la Proposition \ref{condiint}.  La paire $(t,u)$ est uniquement d\'etermin\'ee par ces propri\'{e}t\'{e}s, pour $u$ dans un voisinage de $u_0\equiv 0$.
  \end{itemize} 
  \end{theo}
  \section{Preuve du Th\'eor\`eme \ref{casoliso1}, un argument de transversalit\'e}
  Soit $\{f_s\}_{s\in \Sigma\subset \C }$ une famille lisse \`a un param\`etre complexe de fonctions dans $\gaminf$.  D\'efinissons l'application d'\'evaluation, qui est de classe $\cinf$ mais  pas n\'ecessairement bonne
  \begin{eqnarray*}
   e:\Sigma &\longrightarrow& \gaminf\\
   s&\longmapsto&f_s.
   \end{eqnarray*}
   Si la famille $\{f_s\}$ v\'erifie aussi que $\|f_s-\lambda z\|_r<\bar{\eps}$ pour tout $s$ dans un certain disque $D(0,R) \subset \Sigma \subset \C$, o\`u $\bar{\eps},r$ sont ceux du Th\'eor\`eme \ref{casoliso3}, nous pouvons d\'efinir  une application  $t:D(0,R)\to \C$  de classe $\cinf$, et une application $u:D(0,R)\to \Cinf$ de classe $\cinf$ par  
  \begin{eqnarray*}
  s&\mapsto& t(s)=\pi_1\Xi\big(e(s)\big)\\
  s&\mapsto& u(s)=\pi_2\Xi\big(e(s)\big)
  \end{eqnarray*}
  o\`u l'application $\Xi$ est celle fournie par le Th\'eor\`eme \ref{casoliso3}. Dans ce cas  la courbe $u(s)$ est invariante par la dynamique holomorphe fibr\'ee $F_s(\theta,z)=(\theta+\alpha,e^{t(s)}f_s(\theta,z))$, est de degr\'e nul  et son nombre de rotation transversal est $\varrho_{tr}(u(s))=\beta$. Le but de cette section est de montrer que sous  l'hypoth\`ese de transversalit\'e de la famille $\{f_s\}_{s\in \Sigma}$ on peut choisir un rayon $R>0$ et trouver un param\`etre $s^*$ dans $D(0,R)$ tel que $t(s^*)=0$, c'est \`a dire, qu'on n'a pas besoin de faire une correction $e^t$ sur la dynamique $F_{s^*}(\theta,z)=(\theta+\alpha,f_{s^*}(\theta,z))$ afin d'obtenir la courbe invariante avec le bon nombre transversal de rotation. 
  \paragraph{L'Indice de Kronecker. }
  Nous allons faire dans cette section un petit rappel sur l'indice de Kronecker d'une fonction du plan $\R^2$ vers $\R^2$ par rapport \`a un disque. Soit une fonction $G:D\subset \R^2\to \R^2$ de classe $C^2$, qui s'\'ecrit $G(x)=(g_1(x),g_2(x))$ d\'efinie sur un disque  $D$ du plan. Nous d\'efinissons l'indice de Kronecker de $G$ sur le bord $\partial D$ par l'int\'egrale
  \begin{equation*}
  n(G;D)=\frac{1}{2\pi}\int_{\partial D}\frac{g_1dg_2-g_2dg_1}{g_1^2+g_2^2}
  \end{equation*}   
  quand il n'y a pas de z\'eros de $G$ sur le bord $\partial D$. Ce nombre  mesure le nombre de tours que la courbe $G(\partial D)$ fait au tour de z\'{e}ro. En fait, il n'est pas difficile  de voir que 
  \begin{equation*}
  n(G;D)=\frac{1}{2\pi}\int_{\partial D}G^*(\delta \theta)=\frac{1}{2\pi}\int_{G\circ\partial D}\delta \theta,
  \end{equation*} 
  o\`u $\delta \theta$ est la 1-forme d'\'{e}l\'{e}ment d'angle autour de l'origine complexe. Parmi les diverses propri\'et\'es de l'indice de Kronecker nous allons utiliser seulement celles qui sont contenues dans la proposition suivante, dont la preuve se trouve dans \cite{ELON81}:
  \begin{prop}\label{krone}
  Soit $G:D\subset \R^2\to \R^2$ et $J:D\subset \R^2\to \R^2$ deux fonctions de classe $C^2$. On a
    \begin{enumerate}
  \item Si  $0$ n'appartient pas au segment $[G(z),J(z)]$ pour tout $z$ dans $\partial D$ alors $n(G; D)=n(J; D)$.
  \item Si $G$ n'a pas de z\'{e}ros dans le disque $D$ alors  $n(G;D)=0$.
  \item Si $|G(z)|<|J(z)|$ et $J(z)\neq 0$ pour tout $z$  dans $\partial D$ alors $n(J+G; D)=n(J; D)$.
  \item Si $G$ a un z\'ero unique et non d\'eg\'en\'er\'e dans $D$ alors $n(G; D)=\pm1$.
    \end{enumerate}
  \end{prop}
  On dit qu'un point $z$ dans  $D$ est un z\'ero non d\'eg\'en\'er\'e de $G$ si le d\'eterminant de la matrice Jacobienne $J(G)$ de $G$ est non nul.  Il est clair que $2., 3.$ sont une  cons\'{e}quence de $1.$ Nous allons utiliser cette proposition pour montrer l'existence d'un z\'ero pour une fonction qui n'est pas tr\`{e}s bien comprise, en la comparant \`a une autre qui poss\`{e}de un z\'ero unique et non d\'eg\'en\'er\'e, dont l'indice de Kronecker n'est pas nul. 
  \subsection{La fonction t(s)}
  Nous allons \'etudier par la suite le comportement de la fonction  $t(s)$. Nous savons d'apr\`es (\ref{formulat}) qu'on a explicitement 
  \begin{equation*}
  t(s)=2\pi i\beta-\int_{\T}\log f'_s\big(\theta,u(s)(\theta)\big) d\theta
  \end{equation*}
  d'o\`u on peut calculer la diff\'erentielle de $t$ au point $s$ dans la direction $\Delta s$ par
  \begin{equation*}\label{dts}
  Dt(s)\Delta s=-\int_{\T}\frac{\partial_s f_s'\big(\theta,u(s)(\theta)\big)\Delta s+f_s''\big(\theta,u(s)(\theta)\big)Du(s)\Delta s(\theta)}{f'_s\big(\theta,u(s)(\theta)\big)}d\theta.
  \end{equation*}
  Le Th\'eor\`eme des Fonctions Implicites de Hamilton nous fournit aussi l'expression de la diff\'erentielle de la fonction implicite engendr\'ee, ce qui nous permet de calculer
  \begin{eqnarray*}
  Du(s)\Delta s&=& D\pi_2\Xi\big(e(s)\big)\big[De(s)\Delta s\big]\\
  &=&D\pi_2\Xi(f_s)\big[\partial_sf_s\Delta s\big]\\
  &=&-\Big(D_2\Theta\big(f_s,\pi_2\Xi(f_s)\big)\Big)^{-1}\Big[D_1\Theta\big(f_s,\pi_2\Xi(f_s)\big)[\partial_sf_s\Delta s]\Big].
  \end{eqnarray*}
  Ceci nous permet d'affirmer que s'il existe une constante $T>1$ telle que $\|\partial_sf_s\|_0<T$ pour tout $s$ dans $D(0,R)$ alors il existe une constante positive $T'$, qui d\'epend seulement de la paire $(\alpha,\beta)$ et de $T$, telle que 
  \begin{equation*}
  \|Du(s)\Delta s\|_0<T'|\Delta s|.
  \end{equation*}
  Nous pouvons \'ecrire donc $Dt(s)=w(s)+A(s)$ avec $w(s)$ une matrice $\in \mathbb{M}_2(\R)$ et la matrice $A(s)\in \mathbb{M}_2(\R)$ de taille contr\^ol\'ee par la taille de $f''_s$. Plus pr\'{e}cis\'{e}ment
  \begin{eqnarray*}
  w(s)&=&-\int_{\T}\Big(f'_s\big(\theta,u(s)(\theta)\big)\Big)^{-1} \partial_s f_s'\big(\theta,u(s)(\theta)\big)d\theta\\
  A(s)\Delta s&=&\int_{\T} \Big(f'_s\big(\theta,u(s)(\theta)\big)\Big)^{-1}f''_s\big(\theta,u(s)(\theta)\big)Du(s)\Delta s (\theta)d\theta.
      \end{eqnarray*}
      Si on suppose qu'il existe des nombres r\'eels positifs $\eps_0,\eps_1,\eps_2$ tels que $\|f_s\|_r\leq \eps_0,\|f_s'-\lambda\|_r\leq \eps_1, \|f''_s\|_r\leq \eps_2$ et $\eps_0+\eps_1+\eps_2\leq\bar{\eps}$   alors on a que  
          \begin{equation}\label{AAA}
      \|A(s)\|_{\mathfrak{L}}<2T'\eps_2.
      \end{equation}
       Nous pouvons contr\^oler aussi la distance $\|w(s)-w(0)\|_{\LL}$ sur le bord d'un disque $D(0,R)$. Pour cela nous consid\'erons les in\'{e}galit\'{e}s suivantes:
      \begin{equation*}
      \big\|\partial_sf'_s(\theta,z)-\partial_sf'_s\big|_{s=0}(\theta,z)\big\|_{\LL}\leq \|\partial^2_sf'_s\|_0|s|\leq TR
      \end{equation*}
       si on suppose que $\|\partial^2_sf'_s\|_0\leq T$. De fa\c con similaire nous avons 
      \begin{equation*}
      \big\|\partial_sf'_s\big|_{s=0}(\theta,u(0)(\theta))-\partial_sf'_s\big|_{s=0}(\theta,0)\big\|_{\LL}\leq \|\partial_sf''_s\big|_{s=0}\|_0\|u(s)\|_0\leq TT'R
      \end{equation*}
      si $\|\partial_sf''_s\|_0\leq T$. Il existe donc une constante positive $T''$, qui d\'epend seulement de $T$ et $(\alpha,\beta)$,  telle que
      \begin{equation}\label{WWW}
      \|w(s)-w(0)\|_{\LL}\leq T''R.
      \end{equation}
      \subsection{Transversalit\'e et fin de la preuve}
      Consid\'erons la fonction affine $W(s)=t(0)+w(0)s$. L'hypoth\`ese de transversalit\'{e} 
      \begin{equation*}\label{LLL}
      \bigg [\partial_s\int_{\T}f'_s(\theta,0)d\theta\Big|_{s=0} \bigg]_{\LL}>L^{-1},
      \end{equation*} 
     pour une constante $L>1$, impliquera que $W$ cro\^{\i}t suffisamment vite, ce qui va nous permettre de la comparer \`a la diff\'erence $t(s)-W(s)$ sur le bord d'un disque $D(0,R)$, de rayon $R$ assez grand de fa\c con   que le z\'ero     (unique et non d\'{e}g\'{e}n\'{e}r\'{e}) 
     \begin{equation*}
      \tilde{s}=-w(0)^{-1}t(0) 
      \end{equation*}
      de $W$ ,  soit contenu dans le disque (nous allons utiliser le point $3.$ de la Proposition \ref{krone}). D'apr\`es le Lemme \ref{cinf2lemaA} et l'estimation (\ref{estimacionT}) la taille de $\tilde{s}$ est born\'ee par $LC\eps_1$.  D'apr\`es (\ref{AAA}) et (\ref{WWW}) la diff\'erence entre les fonctions $W$ et $t$ sur le bord du disque $D(0,R)$ est born\'ee par
     \begin{eqnarray*}
     |W(s)-t(s)|&\leq& R\|D(W-t)\|_{\mathfrak{L}}\\
     &\leq& R\big( \|w(s)-w(0)\|_{\LL}+\|A(s)\|_{\mathfrak{L}}\big)\\\label{wmenost}
     &\leq& R(T''R+2T'\eps_2).
     \end{eqnarray*}
      D'autre part la taille de $W$ est minor\'ee sur le bord du disque par
      \begin{equation*}\label{tallaW}
      |W(s)|>RL^{-1}-C\eps_1.
      \end{equation*}
      Pour pouvoir utiliser le point $3.$ de la Proposition \ref{krone} il suffira alors que l'in\'egalit\'e suivante soit v\'erifi\'ee
      \begin{equation*}\label{desigualdadw}
      RL^{-1}-C\eps_1>T''R^2+2T'\eps_2R
      \end{equation*}
    et c'est justement le cas si nous prenons $R=2LC\eps_1$ et $\bar{\eps}$ suffisamment petit. Ce choix nous permet aussi d'assurer la pr\'{e}sence de $\tilde{s}$ dans $D(0,R)$ et de cette mani\`{e}re la pr\'{e}sence d'un z\'ero $s^*$ de $t$ dans $D(0,R)$, car il n'y a que deux possibilit\'{e}s, soit $t$ s'annule sur le bord du disque, ou bien l'indice $n(t,D(0,R))=\pm1$ et le point $2.$ de la Proposition \ref{krone} nous donne l'existence de $s^*$. Nous venons de montrer ainsi la
    \begin{prop}\label{casicasi}
    Pour toute paire $(\alpha,\beta)$ qui v\'erifie l'hypoth\`ese $\C\D_1$ et pour toutes constantes $L>1,T>1$ il existe $\bar{\eps}>0$, un nombre naturel $r\geq 2$, qui d\'epend seulement de la paire $(\alpha,\beta)$, et une constante positive universelle $C$ tels que si une famille lisse \`a un param\`etre complexe $\{f_s\}_{s\in \Sigma\subset \C}$ de fonctions de $\T$ vers  $\gaminf$ v\'erifie  pour  des nombres r\'eels positifs $\eps_0,\eps_1,\eps_2$,  $\eps_0+\eps_1+\eps_2\in (0,\bar{\eps}]$, 
    \begin{itemize}
    \item[$\bullet$] $\bigg[ \partial_s\int_{\T}f'_s(\theta,0)d\theta\Big|_{s=0}\bigg]_{\LL}>L^{-1}$
    \item[$\bullet$] $\|f_s\|_r\leq \eps_0,\|f'_s-\lambda\|_r\leq \eps_1,\|f''_s\|_r\leq \eps_2$ pour tout $s$ dans $D(0,2CL\eps_1)$
    \item[$\bullet$] $\|\partial_s^2f'_s\|_0+\|\partial_sf''_s\|_0\leq T$ pour tout $s$ dans $D(0,2CL\eps_1)$
    \end{itemize}
    alors il existe un param\`{e}tre $s^*$ dans le disque $D(0,2CL\eps_1)$ et une courbe $u:\T\to \D$ de classe $\cinf$  qui est invariante par la dynamique holomorphe fibr\'ee $F_{s^*}(\theta,z)=(\theta+\alpha,f_{s^*}(\theta,z))$, est de degr\'e nul  et son nombre de rotation transversal est $\varrho_{tr}(f_{s^*},u)=\beta$.
    \end{prop}
        Nous pouvons finir \`a pr\'esent  la preuve du Th\'eor\`eme \ref{casoliso1}. Pour cela, supposons qu'il existe $\tilde{\eps}>0$ et  une constante $M>1$ tels que 
    \begin{eqnarray*}
    \|\rho_{0,s}\|_r&\leq &\tilde{\eps}\\
    \|\rho_{1,s}\|_r&\leq&\tilde{\eps}\\
    \|\rho_s\|_r&\leq& M.
    \end{eqnarray*}
     Nous faisons un changement d'\'echelle sur $z$ de taille $m>1$, c'est \`a dire, nous d\'efinissons une nouvelle famille lisse $\{\tilde{f}_s\}_{s\in \Sigma\subset \C}$ par
    \begin{equation*}
    \tilde{f}_s(\theta,z)=mf_s(\theta,m^{-1}z)
    \end{equation*}
    ce qui avec la notation usuelle nous  donne
    \begin{eqnarray*}
    \|\tilde{\rho}_{0,s}\|_r&\leq& m\tilde{\eps}\\
    \|\tilde{\rho}_{1,s}\|_r&\leq& \tilde{\eps}\\
    \|\tilde{\rho}_s\|_r&\leq& m^{-1}M.
    \end{eqnarray*}
    Si nous posons  $m=3M\bar{\eps}^{-1}$ (quitte \`a diminuer $\bar{\eps}$ on a que $m>1$) et $\tilde{\eps}\leq \frac{\bar{\eps}^2}{9M}$ nous avons 
    \begin{eqnarray*}
    \|\tilde{\rho}_{0,s}\|_r&\leq& \frac{\bar{\eps}}{3}\\
    \|\tilde{\rho}_{1,s}\|_r&\leq& \tilde{\eps} < \frac{\bar{\eps}}{3}\\
    \|\tilde{\rho}_s\|_r&\leq& \frac{\bar{\eps}}{3}.
    \end{eqnarray*}
    La Proposition \ref{casicasi} nous donne donc un param\`{e}tre $s^*$ dans le disque $D(0,2CL\tilde{\eps})$ et une courbe $u$ de classe $\cinf$ qui est invariante par la dynamique holomorphe fibr\'ee $\tilde{F}_{s^*}$, ce qui signifie
    \begin{eqnarray*}
    \tilde{F}_{s^*}(\theta,u(\theta))&=&\big(\theta+\alpha,u(\theta+\alpha)\big)\\
    \big(\theta+\alpha,mf_{s^*}(\theta,m^{-1}u(\theta))\big)&=&\big(\theta+\alpha,u(\theta+\alpha)\big)\\
    f_{s^*}\big(\theta,m^{-1}u(\theta)\big)&=&m^{-1}u(\tal).
        \end{eqnarray*}
        On voit ainsi  que la courbe $m^{-1}u$ est invariante par la dynamique holomorphe fibr\'ee $F_{s^*}$. En plus $\tilde{f}'_s(\theta,z)=f'_s(\theta,m^{-1}z)$, donc le nombre de rotation transversal de la courbe $m^{-1}u$ est  aussi \'egal \`a $\beta$ ce qui termine la preuve du Th\'eor\`eme \ref{casoliso1}$\quad_{\blacksquare}$
\bibliographystyle{plain}
\bibliography{cinftesis.bib}
\end{document}